\documentclass[reqno]{amsart}
\usepackage{amsmath,amsfonts,amssymb,amsthm}
\usepackage{eucal}
\usepackage[latin1]{inputenc}
\newcommand{\R}{\Bbb R}

\newcommand{\Chi}{{\bf \chi}}
\newtheorem{Theorem}{Theorem}
\newtheorem{Lemma}{Lemma}
\newtheorem{Remark}{\small Remark}
\newtheorem{Proposition}{Proposition}

\newenvironment{Proof}{\smallskip \noindent{\bf Proof}: \,\,}

\begin{document}
\title[Well-posed]{Well-posedness for a higher order nonlinear Schr\"{o}dinger 
equation in Sobolev spaces of negative indices.}
\author{Xavier Carvajal}
\address{IMECC-UNICAMP, Caixa
Postal:6065, 13083-859, Bar\~ao Geraldo, Campinas, SP,  Brazil}
\email{carvajal@ime.unicamp.br}
\thanks{This research was supported by FAPESP, Brazil.}
\begin{abstract} We prove that, the initial value problem associated to
\begin{equation} 
\partial_{t}u+i\alpha 
\partial^{2}_{x}u+\beta  \partial^{3}_{x}u+i\gamma|u|^{2}u  =  0,
 \quad x,t \in \R, 
\nonumber
\end{equation}
is locally well-posed in $\mathbf{H^{s}}$ for $s>-1/4$. 
\end{abstract}
\maketitle
\noindent

\section{Introduction}

In this work, we study a particular case of the following initial value problem (IVP)
\begin{eqnarray}\label{0y0}
\partial_{t}u+i\alpha
\partial^{2}_{x}u+\beta \partial^{3}_{x}u+F(u) &=&  0, \quad \,\, x,t \in \R,\\
u(x,0) &=& u_0(x)\nonumber
\end{eqnarray}
where $u$ is a complex valued function,  $F(u)=i\gamma|u|^{2}u+\delta
|u|^{2}\partial_{x}u+\epsilon u^{2}\partial_{x}\overline{u}$, $\gamma , \delta ,
\epsilon  \in \mathbb{C}$ and  $\alpha,\beta \in \R$ are constants.

A. Hasegawa and Y. Kodama \cite{[H-K], [Ko]}, proposed (\ref{0y0}) as a model for propagation of pulse in optical fiber. 
We will study the IVP (\ref{0y0}) in 
Sobolev space $\mathbf{H^{s}(\R)}$ under the condition $\delta=\epsilon=0$,
$\beta \neq 0$ (see case
 {\bf iv)} in Teorema \ref{Teo1} below).
Our definition of local well-posedness includes:
existence, uniqueness, persistence and continuous dependence of solution on given data (i.e. continuity of application $u_{0} \mapsto u(t)$ from $X$ to 
$\mathcal{C}([-T,T];X)$).

If  $T < \infty$ we say that the IVP is locally well-posed in  $X$. If some
hypothesis in the definition of local well-posedness fails, we say that the IVP is ill-posed.

Particular cases of (\ref{0y0}) are the followings:

{\bf $\bullet$} Cubic nonlinear  Schr\"odinger equation 
 (NLS), ($\alpha=\mp 1$, $\beta=0$, $\gamma=-1$,
$\delta=\epsilon=0$).
\begin{align}\label{3y0}
iu_{t} \pm u_{xx} + |u|^{2}u =  0, \quad x,t \in \R.
\end{align}
Best known local result for the IVP associated to (\ref{3y0}) is in $\mathbf{H^{s}(\R)}$, $s \geq 0$, obtained by 
Tsutsumi \cite{[Ts]}.

{\bf $\bullet$} Nonlinear Schr\"odinger equation with derivative ($\alpha=-1 $, $\beta=0$,
$\gamma=0$, $\delta=2\epsilon$).
\begin{align}\label{4y0}
iu_{t} + u_{xx} + i\lambda(|u|^{2}u)_{x}=0, \quad x,t\in \R.
\end{align}
Best known result for the IVP associated to  (\ref{4y0}) is in $\mathbf{H^{s}(\R)}$, $s \geq
1/2$, obtained by Takaoka \cite{[T1]}.

{\bf $\bullet$} Complex modified  Korteweg-de Vries (mKdV) equation  ($\alpha=0$, $\beta=1$, $\gamma=0$, $\delta=1$,
$\epsilon=0$).
\begin{align}\label{83y0}
u_{t}+u_{xxx}+|u|^{2}u_{x}=0, \quad x,t \in \R.
\end{align}
If $u$ is real, (\ref{83y0}) is the usual mKdV equation and
Kenig et al. \cite{[KPV1]}, proved the IVP associated to it is locally well-posed in $\mathbf{H^{s}(\R)}$, $s \geq 1/4$.

Laurey \cite{[C1],[C2]} proved that the IVP associated to (\ref{0y0}) is locally well-posed in $\mathbf{H^{s}(\R)}$, $s > 3/4$.

Staffilani \cite{[G]} improved this result by proving the IVP associated to (\ref{0y0}) is locally well-posed in $\mathbf{H^{s}(\R)}$, $s
\geq 1/4$.

When $\alpha, \beta $ are functions of\, $t$, we proved in \cite{C1, CL1} local well-posedness in $\mathbf{H^{s}(\R)}$, $s \geq 1/4$.
Also we studied in \cite{C1, MP} the unique continuation property for the solution of
(\ref{0y0}).
Regarding the ill-posedness of the IVP (\ref{0y0}), we proved in
\cite{CL0} the following theorem.

\begin{Theorem}\label{Teo1} The mapping data-solution $u_{0} \mapsto
u(t)$ for the IVP (\ref{0y0}) is not $\mathcal{C}^{3}$ at origin in the following cases:\\
 {\bf i)}\,\,\, $\beta=0$, $\alpha \neq 0$, $\delta =\epsilon=0$, $\gamma \neq 0$  for $s<0$. \\
{\bf ii)}\,\,\, $\beta=0$, $\alpha \neq 0$, $\delta \neq 0$ or $\epsilon \neq 0$  for $s<1/2$.\\
{\bf iii)}\,\,\, $\beta \neq 0$, $\delta \neq 0$ or $\epsilon \neq 0$  for $s<1/4$.\\
{\bf iv)}\,\,\, $\beta \neq 0$, $\delta = \epsilon=0$, $\gamma \neq 0$  for $s<-1/4$.\\
\end{Theorem} 
In this work, considering the case {\bf iv)} in Theorem \ref{Teo1}, we prove the following theorem. 
\begin{Theorem}\label{T2}
The IVP associated to {\bf iv)},
\begin{align}\label{0y1}
\partial_{t}u+i\alpha
\partial^{2}_{x}u+\beta  \partial^{3}_{x}u+i\gamma|u|^{2}u =  0, \quad x,t \in \R,
\end{align}
is locally well-posed in $\mathbf{H^{s}(\R)}$, $s > -1/4$.
\end{Theorem}
The following trilinear estimate will be fundamental in the proof of Theorem \ref{T2} 
\\
\begin{Theorem}\label{Tt3}
Let $-1/4 < s \leq 0$, $7/12<b<11/12$, then we have
\begin{align}\label{T3}
\|uv\overline{w}\|_{X^{s,b-1}} \leq C
\|u\|_{X^{s,b}}\|v\|_{X^{s,b}}\|w\|_{X^{s,b}},
\end{align}
where
\begin{align*}\|u\|_{X^{s,b}}=\|<\xi>^{s}<\tau-\phi(\xi)>^{b}\hat{u}
\|_{\mathbf{L^{2}_{\xi}L^{2}_{\tau}}},
\end{align*}
$<\xi>=1+|\xi|$,\, $\phi(\xi)=\alpha \xi^{2}+\beta \xi^{3}$.
\end{Theorem}
\begin{Theorem}\label{T4}
The trilinear estimate (\ref{T3}) fails if $s<-1/4$ and $b\in \R$.
\end{Theorem}
\begin{Remark}
1)
As the equation (\ref{0y0}) preserves $\mathbf{L^{2}}$ norm, the Theorem 
\ref{T2} permits to obtain global existence in $\mathbf{L^{2}}$.

2) From Lemma \ref{lpri1} we note that the value  $7/12+$ is the best possible
for the value very near to $\rho=1/4$, in the trilinear estimate (\ref{T3}).

3) The trilinear estimate is valid for all $s>0$, because it follows by combing
the fact that $<\xi>^{s} \leq \quad
<\xi-(\xi_{2}-\xi_{1})>^{s}<\xi_{2}>^{s}<\xi_{1}>^{s}$ and the estimate (\ref{T3})
 for $s=0$.
 
4) We will use the notation $\|u\|_{\{s,b\}}:=\|u\|_{X^{s,b}}$.

5) When  $\alpha =0, \beta=1$, we have  $-3/4+$ bilinear estimate \cite{[KPV8]}, 
\begin{align*}
\|(uv)_{x}\|_{\{-3/4+,-1/2+\}} \leq C
\|u\|_{\{-3/4+,1/2+\}}\|v\|_{\{-3/4+,1/2+\}}.
\end{align*}

Also we have the 1/4 trilinear estimate \cite{[Tao]},
\begin{align*}
\|(uvw)_{x}\|_{\{1/4,-1/2+\}} \leq C
\|u\|_{\{1/4,1/2+\}}\|v\|_{\{1/4,1/2+\}} \|w\|_{\{1/4,1/2+\}}.
\end{align*}
\end{Remark}
\section{Proof of Theorem \ref{T4}.}
As in \cite{[KPV8]} consider the set 
\begin{align*}
B:=\{(\xi,\tau); N\leq \xi \leq N+N^{-1/2}, |\tau-\phi(\xi)| \leq 1\},
\end{align*}
where $\phi(\xi)=\alpha \xi^{2}+\beta \xi^{3}$. We have $|B| \thicksim N^{-1/2}$.  Let us consider $\hat{v}=\Chi_{B}$, it is not difficult to see that $\|v\|_{\{s,b\}}\leq N^{s}|B|^{1/2}$. Moreover
\begin{align*}
\mathcal{F}(|v|^{2}\overline{v}):=\Chi_{B}\ast\Chi_{B}\ast\Chi_{-B}\gtrsim \frac{1}{N}\Chi_{A},
\end{align*}
where $A$ is a rectangle contained in $B$ such that $|A|\thicksim N^{-1/2}$.\\
 Therefore
\begin{align*}
\|\,|v|^{2}\overline{v}\|_{\{s,b-1\}}=&\|<\xi>^{s}<\tau-\phi(\xi)>^{b-1}\mathcal{F}(|v|^{2}\overline{v})\|_{\mathbf{L_{\xi}^{2}L_{\tau}^{2}}}\\
 \gtrsim & N^{s}\frac{1}{N}N^{-1/4}=N^{s-5/4}.
\end{align*}
As a consequence, for large $N$ the trilinear estimate fails if
$3(s-1/4)<s-5/4$, i.e. if $s<-1/4$. \qed
\\
\section{Proof of Theorem \ref{Tt3}.}
To prove Theorem \ref{Tt3}, we need the following results from elementary calculus. \begin{Lemma}\label{lem1}
\begin{enumerate}
\item If $b>1/2$, $a_{1}, a_{2} \in \R$
\begin{align}\label{el1}
\int_{\R}\frac{dx}{<x-a_{1}>^{2b}<x-a_{2}>^{2b}} \thicksim
\frac{1}{<a_{1}-a_{2}>^{2b}}.
\end{align}
\item If\, $0< c_{1}, c_{2}<1$, $c_{1}+c_{2}>1$, $a_{1}\neq a_{2}$,
then
\begin{align}\label{el2}
 \int_{\R}\frac{dx}{|x-a_{1}|^{c_{1}}|x-a_{2}|^{c_{2}}}\lesssim
 \frac{1}{|a_{1}- a_{2}|^{(c_{1}+c_{2}-1)}}.
\end{align}
\item Let $a \in \R$, $c_{1} \leq c_{2}$, then 
\begin{align}\label{el3}
\frac{|x|^{c_{1}}}{<ax>^{c_{2}}} \leq
\frac{C(c_{1},c_{2})}{a^{c_{1}}},
\end{align}
where $C(c_{1},c_{2})$ is a constant independent of $x$.
\item Let $a, \eta \in \R$, $b>1/2$, then
\begin{align}\label{el4}
\int_{\R} \frac{dx}{<a(x^{2}-\eta^{2})>^{2b}} \lesssim
\frac{1}{|a\eta|}.
\end{align}
\end{enumerate}
\end{Lemma}
\vspace{5mm}

Let $f(\xi,\tau)=\,<\xi>^{s}\,<\tau-\xi^{3}>^{b}\hat{u}$,
\,$g(\xi,\tau)=\,<\xi>^{s}\,<\tau-\xi^{3}>^{b}\hat{v}$,
\,$h(\xi,\tau)=<\xi>^{s}$\, $<\tau-\xi^{3}>^{b}\hat{w}$,
$\eta=(\xi,\tau)$, $x=(\xi_{1},\tau_{1})$, $y=(\xi_{2},\tau_{2})$.\\
We have
\begin{align*}
\|uv\overline{w}\|_{\{s,b\!-\!1\}}= & \, \|\int_{\R^{4}} f(\eta\!+\!x\!-\!y)g(y)\overline{h}(x)K(\eta,x,y)dx
dy \|_{\mathbf{L^{2}_{\eta}}}\\
\leq & \,
\|K(\eta,x,y)\|_{\mathbf{L^{\infty}_{\eta}}\mathbf{L^{2}_{x,y}}}
\|f\|_{\mathbf{L^{2}}}\|g\|_{\mathbf{L^{2}}}\|h\|_{\mathbf{L^{2}}},
\end{align*}
where $r(\xi,\tau)=<\!\xi\!>^{2\rho}<\!\!\tau\!-\!\phi(\xi)\!\!>^{2(1-b)}$,
$\rho=-s$ and
\begin{align*}
K(\eta,x,y)=&\frac{<\!\xi\!+\!\xi_{1}\!-\!\xi_{2}\!>^{2\rho}<\!\xi_{2}\!>^{2\rho}<\!\xi_{1}\!>^{2\rho}}
{r(\xi,\tau)<\!\!\tau_{1}\!-\!\phi(\xi_{1})\!\!>^{2b}<\!\!\tau_{2}\!-\!\phi(\xi_{2})\!\!>^{2b}<\!\!\tau\!+
\!\tau_{1}\!-\!\tau_{2}\!-\!\phi(\xi\!+\!\xi_{1}\!-\!\xi_{2})\!\!>^{2b}}.
\end{align*}
Using (\ref{el1}) we obtain
\begin{align*}
I(\xi,\tau):=\|K(\xi,\tau)\|_{\mathbf{L^{2}_{x,y}}}^{2} \thicksim
\frac{1}{r(\xi,\tau)}\int_{\R^{2}} \frac{G_{\rho}(\xi, \xi_{1},
\xi_{2})\,\,d\xi_{1}d\xi_{2}}
{<\!\tau\!-\!\phi(\xi\!+\!\xi_{1}\!-\!\xi_{2})\!-\!\phi(\xi_{2})\!+\!\phi(\xi_{1})\!>^{2b}},
\end{align*}
where
$G_{\rho}(\xi,\xi_{1},\xi_{2}):=<\!\xi\!+\!\xi_{1}\!-\!\xi_{2}\!>^{2\rho}<\!\xi_{1}\!>^{2\rho}<\!\xi_{2}\!>^{2\rho}$.

For clarity in exposition we consider the case $\alpha=0$,
$\beta=1$, i.e. $\phi(\xi)=\xi^{3}$.  With this consideration we have
\begin{align*}
I(\xi,\tau)\!:=&\frac{1}{<\!\xi\!>^{2\rho}<\!\tau\!-\!
\xi^{3}\!\!>^{2(1\!-b)}}\int_{\R^{2}}
\frac{G_{\rho}(\xi,\!-\xi_{1},\!-\xi_{2})d\xi_{1}d\xi_{2}}{<\!\tau\!-
\!\xi^{3}\!+\!g\!>^{2b}},
\end{align*}
where
$g=g(\xi,\xi_{1},\xi_{2})=3(\xi_{1}-\xi_{2})(\xi-\xi_{1})(\xi+\xi_{2})$.

Supposing  $y=\tau- \xi^{3}$, to get Theorem \ref{Tt3} it is enough to prove 
\begin{Lemma}\label{lpri}
Let $0<\rho <1/4$, $7/12 <b< 11/12$. Then
\begin{align*}
I(\xi,y)\!:=\frac{1}{<\!\xi\!>^{2\rho}<\!y\!\!>^{2(1\!-b)}}\int_{\R^{2}}\!
\frac{G_{\rho}(\xi,\!-\xi_{1},\!-\xi_{2})d\xi_{1}d\xi_{2}}{<\!y\!+\!g(\xi,\xi_{1},\xi_{2})\!>^{2b}}
\leq  C(\rho, b)< \infty,
\end{align*}
where $C(\rho, b)$ is a constant independent of $\xi$ and
$y$.
\end{Lemma}
To prove Lemma \ref{lpri} we need to prove the following lemmas. 
\begin{Lemma}\label{lpri1}
Let $\rho<1/4$, then we have 
\begin{align*}
I(0,0)=& \int_{\R^{2}}
\frac{G_{\rho}(0,-\xi_{1},\!-\xi_{2})d\xi_{1}d\xi_{2}}{<g(0,\xi_{1},\xi_{2})>^{2b}}
= \left\{ \begin{array}{ll} C(\rho, b)< \infty, & \textrm{if $\rho +1/3<b$}\\
\infty , & \textrm{if $\rho+1/3 \geq b$},
\end{array} \right.
\end{align*}
where $C(\rho, b)$ is a constant. \end{Lemma}
\begin{Lemma}\label{lpri2}
Let $\rho <1/4$, $b> 7/12$, then 
\begin{align*}
I(\xi,0)=&\frac{1}{<\!\xi\!>^{2\rho}}\int_{\R^{2}}\!
\frac{G_{\rho}(\xi,\!-\xi_{1},\!-\xi_{2})d\xi_{1}d\xi_{2}}{<\!g(\xi,\xi_{1},\xi_{2})\!>^{2b}}
\leq  C(\rho, b),
\end{align*}
where $C(\rho, b)$ is a constant independent of $\xi$.
\end{Lemma}

In the definition of $I(\xi,y)$ if we make the change of variables $
\xi-\xi_{1}:=\xi \xi_{1}$, $\xi+\xi_{2}:=\xi \xi_{2}$ and
$y=\xi^{3}z$, then  $I(\xi,y)$ becomes
 \begin{align}\label{intpri}
I(\xi,z)=p(\xi,z)\int_{\R^{2}}\frac{H_{\rho}(\xi, \xi_{1},
\xi_{2})d\xi_{1}d\xi_{2}}{<\xi^{3}(z+F(\xi_{1},\xi_{2}))>^{2b}},
\end{align}
where \,\,$p(\xi,z)=\xi^{2}/<\xi^{3}z>^{2(1-b)}<\xi>^{2\rho}$,\,
$F(\xi_{1},\xi_{2})=(2-(\xi_{1}+\xi_{2}))\xi_{1}\xi_{2}$\,\, and $$
H_{\rho}(\xi,\xi_{1},\xi_{2})= 
<\xi(1-(\xi_{1}+\xi_{2}))>^{2\rho}<\xi(1-\xi_{1})>^{2\rho}
<\xi(1-\xi_{2})>^{2\rho}.$$

From here onwards we will suppose $z>0$.
\\
{\bf Proof of Lemma \ref{lpri1}.}

By symmetry it is enough to prove that the following integrals 
\begin{align*}
I_{1}(0,0):=\int_{0}^{\infty}\int_{0}^{\infty}\frac{G_{\rho}(0,-\xi_{1},\!-\xi_{2})d\xi_{1}d\xi_{2}}{<g(0,\xi_{1},\xi_{2})>^{2b}}
,
\,\,\,\,I_{2}(0,0):=\int_{0}^{\infty}\int_{0}^{\infty}\frac{G_{\rho}(0,-\xi_{1},\!\xi_{2})d\xi_{1}d\xi_{2}}{<g(0,\xi_{1},-\xi_{2})>^{2b}}
\end{align*}
are finite.
We will prove  $I_{1}(0,0)$ is finite, the same proof works to prove 
$I_{2}(0,0)$ is finite. Also, by symmetry we can suppose that $0\leq \xi_{2}
\leq \xi_{1}$.

 We have 
\begin{align*}
\int_{1}^{\infty}d\xi_{1}\int_{0}^{\xi_{1}}d\xi_{2}\frac{G_{\rho}(0,-\xi_{1},\!-\xi_{2})}
{<g(0,\xi_{1},\xi_{2})>^{2b}}=&\int_{1}^{\infty}d\xi_{1}\int_{0}^{\xi_{1}/2}d\xi_{2}+
\int_{1}^{\infty}d\xi_{1}\int_{\xi_{1}/2}^{\xi_{1}}d\xi_{2}\\
= I_{1,1}+I_{1,2}.
\end{align*}
As $0\leq \xi_{2} \leq \xi_{1}$, we have $G_{\rho}(0,-\xi_{1},\!-\xi_{2}) \leq
<\xi_{1}>^{4\rho}<\xi_{2}>^{2\rho}$. In $I_{1,1}$ we have $\xi_{1}/2< \xi_{1}-\xi_{2}<\xi_{1}$, therefore
if $b>\rho+1/3$,
\begin{align}
I_{1,1} & \lesssim
\int_{1}^{\infty}<\xi_{1}>^{4\rho}d\xi_{1}\int_{0}^{\xi_{1}/2}\frac{<\xi_{2}\!>^{2\rho}d\xi_{2}}
{<\xi_{1}^{2}\xi_{2}\!>^{2b}}\\
\lesssim &
\int_{1}^{\infty}<\xi_{1}>^{4\rho}\Big(\frac{1}{\xi_{1}^{2}}+\frac{1}{\xi_{1}^{2+4\rho}}+
\frac{1}{\xi_{1}^{2+4\rho}}
\int_{1}^{3\xi_{1}^{3}/2}\frac{x^{2\rho}dx}{(1+x)^{2b}}\Big)d\xi_{1}
\\
=& C(\rho,b)< \infty.
\end{align}
Analogously we can prove that $I_{1,1}=\infty$ if
$b \leq \rho+1/3$.

In $I_{1,2}$ we have $\xi_{1}/2\leq \xi_{2} \leq\xi_{1}$, so
\begin{align*}
I_{1,2} \lesssim &
\int_{1}^{\infty}<\xi_{1}>^{4\rho}d\xi_{1}\int_{\xi_{1}/2}^{\xi_{1}}\frac{<\xi_{1}-\xi_{2}\!>^{2\rho}d\xi_{2}}
{<(\xi_{1}-\xi_{2})\xi_{1}^{2}\!>^{2b}}\\=&\int_{1}^{\infty}<\xi_{1}>^{4\rho}d\xi_{1}\int_{0}^{\xi_{1}/2}
\frac{<x>^{2\rho}dx} {<\xi_{1}^{2}x\!>^{2b}}\\
=& C(\rho,b), \,\,\,\,\, b>\rho+1/3.
\end{align*}\qed

To prove Lemmas \ref{lpri} and \ref{lpri2}, the following propositions will be useful. 
\begin{Proposition}\label{prop2}
Let $\rho\geq 0$, $b>1/3+2\rho/3$, then we have \begin{align*}
J_{2}=\xi^{2+4\rho}\int_{\R^{2}}\frac{d\xi_{1}d\xi_{2}}{<\xi^{3}(z+F)>^{2b}}
\leq C,
\end{align*}
where  $C$ is a constant independent of $\xi$.
\end{Proposition}
\Proof If $\xi_{1}\leq 0$, $\xi_{2}\leq 0$, then $|z+F|\geq
|\xi_{1}+\xi_{2}||\xi_{1}\xi_{2}|$. Therefore by Lemma \ref{lpri1} and by symmetry, it is enough to consider 
$\xi_{1} \geq 0$. We have
$|z\!+\!F|\!=\!|\xi_{1}||(\xi_{2}\!+\!(\xi_{1}-2)/2)^{2}-\!
(\xi_{1}\!-\!2)^{2}/4-z/\xi_{1}|$. Let $l^{2}=(\xi_{1}\!-\!2)^{2}/4+z/\xi_{1}$, $c(\rho)=(2+4\rho)/3$,
then making change of variable $\eta=\xi_{2}\!+\!(\xi_{1}-2)/2$
and using (\ref{el2}) and (\ref{el3}) we have \begin{align*}
J_{2} =& \xi^{2+4\rho}\int_{0}^{\infty}d\xi_{1}
\int_{\R}\frac{d\eta}{<\xi^{3}\xi_{1}(\eta^{2}-l^{2})>^{2b}}\\
\lesssim &\int_{0}^{\infty}d\xi_{1}\int_{\R}\frac{l \,dx}{[|\xi_{1}|l^{2}|x^{2}-1|]^{c(\rho)}}\\
\lesssim &\int_{0}^{\infty}
\frac{d\xi_{1}}{|\xi_{1}|^{c(\rho)}|\xi_{1}-2|^{(1+8\rho)/3}}\int_{\R}
\frac{dx}{|x^{2}-1|^{c(\rho)}}\\
\lesssim & \,C.
\end{align*}\qed

\begin{Proposition}\label{prop1}
Let $|\xi|>1$, $\rho<1/4$, then 
\begin{align*}
J_{1}=\xi^{2+4\rho}
\int_{0}^{\infty}\xi_{1}^{4\rho}\int_{\R}\frac{d\xi_{1}d\xi_{2}}{<\xi^{3}(z+F)>^{2b}}
\leq C,
\end{align*}
where $C$ is a constant independent of $\xi$.
\end{Proposition}
\Proof  By Proposition \ref{prop2} we can suppose $\xi_{1}>4$, so
$(\xi_{1}-2)> \xi_{1}/2$. Using (\ref{el4}) and making change of variables as above, we have\begin{align*}
J_{1}\lesssim \,
\frac{\xi^{2+4\rho}}{|\xi|^{3}}\int_{4}^{\infty}\frac{\xi_{1}^{4\rho}}{\xi_{1}\,l}
d\xi_{1}
\leq C.
\end{align*}\qed
\\
{\bf Proof of Lemma \ref{lpri2}.}\\
{\bf a) If $|\xi| \leq 1$}.\\
Let\, $A_{1}\!=\!\{(\xi_{1},\xi_{2})/ |\xi_{1}|\!>\!2, |
\xi_{2}|\!>\!2\}$, $A_{2}\!=\!\{(\xi_{1},\xi_{2})/
|\xi_{1}|\!\leq\! 2, | \xi_{2}|\!\leq 2\}$,
$A_{3}=\{(\xi_{1},\xi_{2})/ |\xi_{1}|\leq 2, | \xi_{2}|>2\}$ and
$A_{4}=\{(\xi_{1},\xi_{2})/ |\xi_{1}|>2, | \xi_{2}|\leq 2\}$  and
consider $I(\xi,0)=\sum_{j=1}^{4}I_{j}(\xi,0)$, where 
$I_{j}(\xi,0)$ is defined in the region  $A_{j}$. Obviously 
$I_{2} \leq C$. In $A_{1}$ we have $|\xi-\xi_{1}|>|\xi_{1}|/2$ and
$|\xi+\xi_{2}|>|\xi_{2}|/2$, therefore Lemma \ref{lpri1} gives 
$I_{1} \leq C$. In $A_{3}$ we have  $|\xi+\xi_{2}|>|\xi_{2}|/2$, 
and consequently 
\begin{align*}
I_{3}(\xi,0)\lesssim & \frac{1}{<\xi>^{2\rho}}
\int_{A_{3}}\frac{<\xi_{2}>^{4\rho}d\xi_{1}d\xi_{2}}{<(\xi_{1}-\xi_{2})\xi_{2}(\xi-\xi_{1})>^{2b}}\\
=& \frac{1}{<\xi>^{2\rho}}\int_{A_{3}\cap
\{|\xi_{1}-\xi_{2}|>|\xi_{2}|\}}+
\frac{1}{<\xi>^{2\rho}}\int_{A_{3}\cap \{|\xi_{1}-\xi_{2}|\leq
|\xi_{2}|\}}\\
=&I_{3,1}(\xi,0)+I_{3,2}(\xi,0).
\end{align*}
In the first integral, for $\rho <1/4$, $b>1/2$ we have 
\begin{align*}
I_{3,1}(\xi,0)\lesssim &\frac{1}{<\xi>^{2\rho}}\int_{|\xi_{2}|>2}<\xi_{2}>^{4\rho}d\xi_{2}\int_{|\xi_{1}|\leq2}\frac{d\xi_{1}}{<\xi_{2}^{2}(\xi-\xi_{1})>^{2b}}\\
\lesssim &\frac{1}{<\xi>^{2\rho}}\int_{|\xi_{2}|>2}\frac{<\xi_{2}>^{4\rho}d\xi_{2}}{\xi_{2}^{2}}\\
\leq &C.
\end{align*}
To estimate $I_{3,2}(\xi,0)$ we make the change of variables $\eta_{2}=\xi_{1}-\xi_{2}$, $\eta_{1}=\xi_{1}$ and as $|\xi_{1}|\leq 2$ we obtain the same estimate as that for $I_{3,1}(\xi,0)$.

By symmetry we can estimate $I_{4}$ in the same manner as $I_{3}$.
\\{\bf b) If $|\xi| >1$}.\\
Let us consider $I(\xi,0)$ in the form (\ref{intpri}) and let $B_{1}=\{
|\xi_{1}+ \xi_{2}|>4\}$ and $B_{2}=\{ |\xi_{1}+ \xi_{2}|\leq 4\}$, then $I(\xi,0)=I_{1}(\xi)+I_{2}(\xi)$, where $I_{j}(\xi)$ is defined in $B_{j}$. In $B_{1}$ we have 
\begin{align}\label{B1}
|2-(\xi_{1}+\xi_{2})|>|\xi_{1}+\xi_{2}|/2,\,\,
|1-(\xi_{1}+\xi_{2})|\leq 5|\xi_{1}+\xi_{2}|/4,
\end{align} moreover $B_{1} \subset \{|\xi_{1}| \geq 2\}\cup\{|\xi_{2}| \geq
2\}=:B_{1,1}\cup B_{1,2}$ and therefore $I_{1}(\xi,0)\leq
I_{1,1}(\xi)+I_{1,2}(\xi)$, where $I_{1,j}(\xi)$ is defined in $B_{1,j}\cap B_{1}$. In $B_{1,1}$ we have $|\xi_{1}|/2 \leq
|1-\xi_{1}|\leq 3|\xi_{1}|/2$, therefore using (\ref{B1}), we obtain that $I_{1,1}(\xi) \lesssim I(0,0) \leq C$ if $\rho <1/4$,
$\rho+1/3 <b$. In similar manner we have $I_{1,2}(\xi)\lesssim
I(0,0) \leq C$.

From definition of $B_{2}$ we have $H_{\rho} \lesssim <\xi>^{2\rho}\,
<|\xi|+|\xi||\xi_{2}|>^{4\rho}$, so using symmetry and Propositions \ref{prop2} and \ref{prop1}, we have $I_{2}(\xi) \leq C < \infty$\, if \,$0 \leq \rho <1/4$, $ b>\rho +1/3$.
\\
\\
{\bf Proof of Lemma \ref{lpri}.}
\\
Let $0 \leq \rho <1/4$, $7/12 < b < 11/12$. Using symmetry and Lemma \ref{lpri2} it is enough to prove 
\begin{align*}
J=p(\xi,z)\int_{0}^{\infty}\int_{\R}\frac{H_{\rho}(\xi, \xi_{1},
\xi_{2})d\xi_{1}d\xi_{2}}{<\xi^{3}(z+F(\xi_{1},\xi_{2})>^{2b}}\leq
C < \infty.
\end{align*}
By Lemma \ref{lpri2} we can suppose $|\xi|^{3}z \geq 1$, because if \,$|\xi|^{3}z
< 1$ then $$<\xi^{3}(z+F)>^{-2b} \leq
2^{2b}<\xi^{3}F>^{-2b}.$$ Also  by symmetry  we can suppose $|\xi_{2}|
\leq |\xi_{1}|$.\\ Therefore\begin{align}\label{H}
H_{\rho}(\xi,\xi_{1},\xi_{2})\lesssim
1+|\xi|^{6\rho}+|\xi|^{6\rho}|\xi_{1}|^{6\rho}.
\end{align}
Using Proposition \ref{prop2} we can suppose $|\xi_{1}|>4$ ($l^{-1} \leq
|\xi_{1}|^{-1}$). \\ \\
{\bf a) \,\,If $|\xi| |\xi_{1}| \leq 1$}.\\
We have $ H_{\rho}\lesssim <\xi>^{6\rho}$ and therefore $J \leq
C<\infty$, by Proposition \ref{prop2}.
\\
{\bf b) \,\,If $|\xi| |\xi_{1}| > 1$}.\\
${\bf i)}$ If $|\xi_{1}|^{3} \leq z$, $|\xi_{1}| \leq z^{1/3}$,
we have $ H_{\rho}(\xi,\xi_{1},\xi_{2})\lesssim
1+|\xi|^{6\rho}+|z|^{2\rho/3}|\xi|^{6\rho}|\xi_{1}|^{4\rho}$. Therefore using
(\ref{el4}), in this region we have
\begin{align*}
\frac{\xi^{2+6\rho}|z|^{2\rho/3}}{<\xi^{3}z>^{2(1-b)}}\int_{1/|\xi|}^{|z|^{1/3}}
|\xi_{1}|^{4\rho}d\xi_{1}
\int_{\R}\frac{d\eta}{<\xi^{3}\xi_{1}(\eta^{2}-l^{2})>^{2b}}
\lesssim & \,\frac{\xi^{2+6\rho}|z|^{2\rho/3}}{<\xi^{3}z>^{2(1-b)}|\xi|^{3}}\int_{1/|\xi|}^{\infty}\frac{|\xi_{1}|^{4\rho}d\xi_{1}}{|\xi_{1}|^{2}}\\
\lesssim &
\frac{(|\xi|^{3}z)^{2\rho/3}}{<\xi^{3}z>^{2(1-b)}}\\
\leq &\, C.
\end{align*}
${\bf ii)}$ If $|\xi_{1}|^{3} \geq z$, $|\xi_{1}| \geq z^{1/3}$, we can proceed as follows.\\ 
By Lemma \ref{lpri2} we can suppose  $|z+F| \leq |F|/2$, so $|F|
\leq 2z$, $|(2-(\xi_{1}+\xi_{2}))\xi_{1}\xi_{2}| \leq 2z$.
This implies that $ |1-\xi_{2}||1-(\xi_{1}+\xi_{2})|\lesssim
1+|\xi_{1}|+z^{2/3}$.\\ Therefore \begin{align*}
H_{\rho} \lesssim &\,\,(<\!\xi>^{4\rho}+\!|\xi|^{6\rho})+
|\xi|^{4\rho}|\xi_{1}|^{4\rho}\!+|\xi|^{6\rho}|\xi_{1}|^{2\rho}\!+\!|\xi|^{4\rho}|\xi_{1}|^{2\rho}
 +|\xi|^{6\rho}|\xi_{1}|^{4\rho}
 \\
 &+|\xi|^{4\rho}z^{4\rho/3}+|\xi|^{6\rho}z^{4\rho/3}
 +|\xi|^{6\rho}z^{4\rho/3}|\xi_{1}|^{2\rho}= \sum_{j=1}^{8}l_{j}.
\end{align*}
We have,
\begin{align}\label{4o}
\frac{|\xi|^{6\rho}}{<\xi>^{2\rho}} \leq |\xi|^{4\rho}.
\end{align}

To estimate the term that contains $l_{1}=<\xi>^{4\rho}+|\xi|^{6\rho}$, we use (\ref{4o}) and Proposition \ref{prop2}.

For terms $l_{j}$, $j=2, \ldots,5$,  we use (\ref{4o}) and Propositions \ref{prop2} and \ref{prop1}  if $|\xi|>1$. If $|\xi| <1$, we integrate in the region $\xi_{1}> 1/|\xi|$  as above. 

 In $l_{6}=|\xi|^{4\rho}z^{4\rho/3}$, we have 
\begin{align*}
\frac{|\xi|^{2}|\xi|^{4\rho}z^{4\rho/3}}{<\xi^{3}z>^{2(1-b)}|\xi|^{3}<\xi>^{2\rho}}\int_{z^{1/3}}^{\infty}
\frac{d\xi_{1}}{\xi_{1}^{2}} \lesssim
\frac{1}{(|\xi|^{3}z)^{(1-4\rho)/3}}\leq C.
\end{align*}

We estimate $l_{7}=|\xi|^{6\rho}z^{4\rho/3}$, as in $l_{6}$ using
(\ref{4o}).

Finally in $l_{8}=|\xi|^{6\rho}z^{4\rho/3}|\xi_{1}|^{2\rho}$,
we have \begin{align*}
\frac{|\xi|^{2+6\rho}z^{4\rho/3}}{<\xi>^{2\rho}|\xi|^{3}}\int_{z^{1/3}}^{\infty}
\frac{|\xi_{1}|^{2\rho}d\xi_{1}}{\xi_{1}^{2}} \lesssim
\frac{(|\xi|^{3}z)^{(6\rho-1)/3}}{<\xi^{3}z>^{2(1-b)}} \leq C.
\end{align*}
\qed \vspace{5mm}

\section{Proof of Theorem \ref{T2}.}

Consider a cut-off function $\psi \in \mathcal{C^{\infty}}$,
such that $0\leq \psi \leq 1$,
\begin{align}
\psi(t)=\left\{ \begin{array}{ll}
1& \textrm{if \,$|t| \leq 1$}\\
 0 & \textrm{if \,$|t| \geq 2$},
\end{array} \right.
\end{align}
and let $\psi_{T}(t):=\psi(t/T)$.
To prove Theorem \ref{T2} we need the following result.

\begin{Proposition}\label{prop3}
Let $-1/2<b' \leq 0\leq b \leq b'+1$, $T \in [0,1]$, then \begin{align}\label{eq1}
\|\psi_{1}(t)U(t)u_{0}\|_{\{s,b\}} = & C \|u_{0}\|_{\mathbf{H^{s}}
}\\
\| \psi_{T}(t)\int_{0}^{t}U(t-t')F(t',\cdot))dt'\|_{\{s,b\}} \leq
& C T^{1-b+b'}\|F(u)\|_{\{s,b'\}},\label{eq2}
\end{align}
where $F(u):= i\gamma|u|^{2}u$.
\end{Proposition}
\Proof The proof of (\ref{eq1}) is obvious. The proof of (\ref{eq2}) is practically done in
\cite{[GTV]}. \qed
\\

Let us consider (\ref{0y1}) in its equivalent integral form  \begin{align}\label{int1}
u(t)=U(t)u_{0}- \int_{0}^{t}U(t-t')F(u)(t',\cdot))dt'.
\end{align}
Note that, if for all $t\in \R$, $u(t)$ satisfies:
\begin{align}\label{int2}
u(t)=\psi_{1}(t)U(t)u_{0}- \psi_{T}(t)
\int_{0}^{t}U(t-t')F(u)(t',\cdot))dt',
\end{align}
then $u(t)$ satisfies (\ref{int1}) in $[-T,T]$. Let $a>0$ and
\begin{align}
X_{a}=\{v \in X^{s,b}; \|v\|_{s,b}\leq a\}.
\end{align}
For $v \in X_{a}$  fixed, let us  define 
\begin{align*}
\Phi(v)=\psi_{1}(t)U(t)u_{0}- \psi_{T}(t)
\int_{0}^{t}U(t-t')F(v)(t',\cdot))dt'.
\end{align*}
Let $\epsilon=1-b+b'>0$, using Proposition \ref{prop3} and  Theorem \ref{Tt3} we obtain 
\begin{align*}
\|\Phi(v)\|_{s,b} \leq &C \|u_{0}\|_{\mathbf{H^{s}}}+CT^{\epsilon}
\|F(v)\|_{s,b'}\\
\leq &\,C\|u_{0}\|_{\mathbf{H^{s}}}+CT^{\epsilon}M^{3}\\
 \leq &\, M,
\end{align*}
where we took $M=2C\|u_{0}\|_{\mathbf{H^{s}}}$, $T^{\epsilon} \leq
1/(2CM^{2})$.

We can prove that $\Phi$ is a contraction in an analogous manner. The proof of
the Theorem ~1
follows by using a standard argument, see for example \cite{[KPV1], [KPV8]}.  

\end{document}